# The Lattice of Closure Relations on a Poset

Michael Hawrylycz
Department of Mathematics
Massachusetts Institute of Technology
Cambridge, MA 02139

and

Victor Reiner
Department of Mathematics
University of Minnesota
Minneapolis, MN 55455

**Abstract**

In this paper we show that the set of closure relations on a finite poset $P$ forms a supersolvable lattice, as suggested by Rota. Furthermore this lattice is dually isomorphic to the lattice of closed sets in a convex geometry (in the sense of Edelman and Jamison [EJ]). We also characterize the modular elements of this lattice (when $P$ has a greatest element) and compute its characteristic polynomial.

## 1 The lattice of closure relations

Let $P$ be a poset. A *closure relation* on $P$ is a map $H : P \to P$ such that for all $x, y \in P$ :

1. $x \leq H(x)$

2. $x \leq y \Rightarrow H(x) \leq H(y)$ (monotone)

3. $H(H(x)) = H(x)$ (idempotent)

If $P = 2^S$, the set of all subsets of a set S, then we say $H$ is a closure on the set S. The structure of all closures on a set has been extensively studied by Ore [Or1,Or2].

We may partially order the set of all closure relations on a poset $P$ by setting $H \leq K$ if $H(x) \leq K(x)$ for all $x \in P$. A closure relation $H$ on $P$ may be regarded as a partition of the elements of $P$, by considering $x, y \in P$ to be in the same block of $H$ if $H(x) = H(y)$. If $H(x) \leq K(x)$ for all $x \in P$, then since $x \leq H(x) \leq K(x)$ for all $x \in P$, $H(x)$ must be in the same *block* as $x$ of $K$ (regarding $K$ as a partition of $P$.) Hence $H$ is a refinement of $K$ as partitions. The converse is evident, and so we have the useful observation that $H \leq K$ if and only if $H$ is a refinement of $K$ as paritions.

For the sake of simplifying the proofs, from now on, we will assume that the poset $P$ is *finite*. However, many of our results are valid under other finiteness conditions on $P$ (such as when $P$ has no infinite chains or no infinite ascending chains).

**Theorem 1.** The partial order of all closure relations on any poset $P$ forms a lattice, denoted $LC(P)$. It is a *join-sublattice* of the lattice $\Pi(P)$ of partitions of the elements of $P$, i.e the join in $LC(P)$ is the same as the join in $\Pi(P)$.

*Proof.* Let $\{H_\alpha\}$ be a set of closure relations on $P$. Then each $H_\alpha$ represents a partition of the elements of $P$. Let $\bigvee H_\alpha$ denote the join of the $\{H_\alpha\}$ in the partition lattice $\Pi(P)$. We will first show that each block of $\bigvee H_\alpha$ has a unique greatest element. Then we can define a map $J : P \to P$ by sending every $x \in P$ to the greatest element in its block of $\bigvee H_\alpha$. Finally we show that $J$ obeys axiom 2 and so is a closure relation, and hence $J$ must be the join of the set $\{H_\alpha\}$ in $LC(P)$.

To show each block of $\bigvee H_\alpha$ has a greatest element, suppose we have two elements $a, b$ that are maximal in some block. By definition of join in $\Pi(P)$, this means that there exists a sequence $B_1, B_2, \ldots, B_n$ of blocks from the various $H_\alpha$'s, such that $a \in B_1, b \in B_n$, and $B_k \cap B_{k+1} \neq \phi$. Let $z_k$ be the maximum element of $B_k$. We have $z_1 = a$ and $z_n = b$ by maximality, so if we prove that $z_k \leq a$ for all $k$, then we will have the contradiction $b \leq a$. Assume by induction that $z_k \leq a$ (true for $k = 1$), let $w \in B_k \cap B_{k+1}$, and assume $B_{k+1}$ is a block of the closure relation $H_i$. Then

$$z_{k+1} = H_i(w) \leq H_i(z_k) \leq H_i(a) = a.$$

since if $H_i(a) > a$ then $a$ would not be maximal.

Now we show $J$ (as defined above) is a closure relation on $P$. Since $J$ sends each element to the maximum element in its block in $\bigvee H_\alpha$, $x \leq J(x)$ and $J(J(x)) = J(x)$. Suppose $x \leq y$ in P. Define a sequence $x_1, x_2, \ldots \in P$ as follows: Let $x_0 = x$, and define $x_{k+1}$ by choosing a closure $H_k$ among $\{H_\alpha\}$ satisfying $H_k(x_k) \neq x_k$, and then setting $x_{k+1} = H_k(x_k)$. By finiteness, this process must stop at some $x_n$, and $x_n = J(x)$. Define the sequence $y_1, y_2, \ldots \in P$ by applying the same sequence of $H_k$'s to $y$, and note that $x_k \leq y_k$ for all $k$ since each $H_k$ is a closure. Thus we have $J(x) = x_n \leq y_n \leq J(y)$. Thus $J$ is monotone, and hence defines a closure relation.

So far we have shown that $LC(P)$ is a join-sub-semilattice of $\Pi(P)$. To show it is a lattice, it suffices to note that it has a minimum element, viz. the *identity closure*, $I(x) = x \ \forall x \in P$, which must be less than any closure on $P$. $\square$

The following proposition characterizes when $LC(P)$ is a sublattice of $\Pi(P)$ for posets $P$ with a greatest element.

**Proposition 2.** Let $P$ have a greatest element $\hat{1}$. Then $LC(P)$ is a sublattice of $\Pi(P)$ if and only $\hat{0} + P$ is a lattice, where $\hat{0} + P$ is the poset obtained by adjoining a new least element $\hat{0}$ to $P$.

Remark: A similar (but harder to state) proposition holds even if $P$ has no greatest element.

*Proof.* ($\Rightarrow$): If $\hat{0} + P$ is not a lattice then there must exist four elements $a, b, c, d \in P$ with $c, d$ both maximal lower bounds for $a$ and $b$. Define the closure $H_a$ by $H_a(x) = a$ if $x \leq a$ and $H_a(x) = \hat{1}$ otherwise. Define $H_b$ similarly. If we form $H_a \wedge H_b$ in $\Pi(P)$ then both $c$ and $d$ will be maximal in the same block. Hence this partition does not correspond to a closure relation. Thus $H_a$ and $H_b$ cannot have the same meet in both $\Pi(P)$ and $LC(P)$.

($\Leftarrow$) Suppose $\hat{0} + P$ is a lattice and let $\{H_\alpha\}$ be a set of closure relations on $P$. Since the meet operation in $LC(P)$ is precisely intersection of the blocks of the individual closure relations, their meet in $\Pi(P)$ and $LC(P)$ coincide if and only if every block in their meet in $\Pi(P)$ has a greatest element. So suppose this is not the case, i.e. let $a, b$ be maximal in some block and $H_\alpha(a) = H_\alpha(b) = c_\alpha, \forall \alpha$. Let $d$ be the greatest lower bound of the $c_\alpha$ in the lattice $\hat{0} + P$. Then $c_\alpha \geq d \geq a, b$ for all $\alpha$, which implies that $H_\alpha(d) = c_\alpha$ for all $\alpha$. Hence $d$ is in the same block as $a$ and $b$, contradicting their maximality. □

## 2 Convexity and mlb-closure

In this section we relate the lattice of closure relations to the notion of a convex geometry as studied by Edelman and Jamison [EJ]. The following proposition is the crucial observation necessary for what follows:

**Proposition 3.** Let $H$ be a closure relation on $P$, and $A$ its set of closed elements, i.e.

$$A = \{x \in P : H(x) = x\}.$$

Then for any subset $B \subseteq A$, all maximal lower bounds of $B$ are in $A$ (we call such a set $A$ *mlb-closed*). Conversely, any $A \subseteq P$ which is mlb-closed defines the closed elements of a unique closure relation $H$.

*Remark:* Note that taking $B$ to be the empty set, we have that any mlb-closed set $A$ must contain all of the maximal elements of $P$.

*Proof.* Let $H$ be a closure relation on $P$, with closed elements $A$. If $B \subseteq A$ has some maximal lower bound $x$, then
$$x \leq H(x) \leq H(b) = b$$
for all $b \in B$. Hence $H(x)$ is also a lower bound of $B$, so by maximality, $H(x) = x$ and $x \in A$. Conversely, given a set $A$ which is mlb-closed, we claim that for any $x$ in $P$ there exists a unique least element of $A_{\geq x}$, where $A_{\geq x}$ is defined to be $\{a \in A : a \geq x\}$. To prove this claim, note that since $A$ contains the maximal elements of $P$, $A_{\geq x} \neq \phi$. If there were two such minimal elements $y, z \in A_{\geq x}$, then they would have a maximal lower bound $w$ above $x$, contradicting their minimality. Thus, given any mlb-closed set $A$, we define $H(x)$ to be the minimum of $A_{\geq x}$. It remains to show that $H$ is the unique closure relation with closed elements $A$. That $H$ satisfies axioms 1 and 3 is evident. To show axiom 2: if $x \leq y$, then we have $A_{\geq x} \supseteq A_{\geq y}$ and hence
$$H(x) = minA_{\geq x} \leq minA_{\geq y} = H(y).$$

Finally, any closure relation is defined uniquely by specifying its closed elements so $H$ is unique. □

Denote by $\overline{A}$ the mlb-closure of a set $A \subseteq P$, i.e. $\overline{A}$ is the intersection of all mlb-closed sets which contain $A$. It can easily be shown that this defines a closure on the set $P$ (i.e. a closure relation on the poset $2^P$). Let $L_{mlb}$ denote the lattice of mlb-closed subsets of $P$ ordered under inclusion, and let $M \subseteq P$ be the set of maximal elements of $P$. We then have

**Theorem 4.** The order-dual of $LC(P)$ is isomorphic to the interval between $\overline{M}$ and $P$ in $L_{mlb}(P)$.

*Proof.* Define a map $f : LC(P) \to 2^P$ by
$$f(H) = \{x \in P : H(x) = x\}.$$

Figure 1: A six-element poset $P$ and its lattice $LC(P)$.

Proposition 3 shows that the image of $f$ is exactly the interval from $\overline{M}$ to $P$ in $L_{mlb}(P)$. It remains then to show that $H \leq K$ if and only if

$$\{x \in P : K(x) = x\} \subseteq \{x \in P : H(x) = x\}.$$

Clearly, if $H \leq K$ and $K(x) = x$ then

$$x = K(x) \geq H(x) \geq x$$

so $H(x) = x$. Conversely, if

$$\{x \in P : K(x) = x\} \subseteq \{x \in P : H(x) = x\},$$

then by the proof of Proposition 3, we have

$$H(y) = min(\{x \in P : H(x) = x\}_{\geq y}) \leq min(\{x \in P : K(x) = x\}_{\geq y}) = K(y).\square$$

**Example.** Consider the poset $P$ shown in figure 1. Here $M = \{a, b\}$ and $\overline{M} = \{a, b, c\}$. The lattice $LC(P)$ is shown with the closure relations pictured as partitions of the set $P$.

In [EJ], Edelman and Jamison define a closure $\bar{\phantom{x}}$ on a set $S$ to be a *convex closure* if it satisfies the following *anti-exchange axiom*: Given distinct $x, y \in S$ and a closed set $A = \overline{A} \subset S$ with $x, y \notin A$, we have

$$x \in \overline{y \cup A} \Rightarrow y \notin \overline{x \cup A}.$$

**Proposition 5** Mlb-closure is a convex closure on the set $P$.

*Proof.* As before let $\overline{A}$ denote the mlb-closure of $A$, and let $x$ and $y$ be distinct elements not in $A$, a closed set. If $x \in \overline{y \cup A}$ then $x$ is a maximal lower bound of some set of maximal lower bounds of some set of maximal lower bounds, etc., of some subsets of $y \cup A$. In this expression for $x$, if $y$ never appears then $x \in \overline{A} = A$, a contradiction. Since $y$ does appear, $x \leq y$. Thus we cannot have $y \leq x$, and so $y \notin \overline{x \cup A}.\square$

Edelman and Jamison [EJ] show that convex closures enjoy a number of interesting properties. Corollary 6 states a few properties of $LC(P)$ which are consequences of mlb being a convex closure. The reader is referred to [EJ] for proofs.

**Corollary 6**

1. $LC(P)$ is a *join-distributive* lattice, i.e. every atomic interval in $LC(P)$ is a Boolean algebra.

2. $LC(P)$ is upper-semimodular, and consequently ranked.

3. If $K \in LC(P)$ is defined by the mlb-closed set $A \subseteq P$, then the rank function of $LC(P)$ is given by $r(K) = card(P) - card(A)$, where $card(A)$ is also the number of blocks in $K$ when regarded as a partition of the elements in $P$.□

We note one further consequence not given in [EJ]. Recall the definition of the *characteristic polynomial* of a ranked finite poset $Q$ with $\hat{0}, \hat{1}$ and rank function $r$:

$$\chi(Q;\lambda) = \sum_{q \in Q} \mu(\hat{0},q)\lambda^{r(\hat{1})-r(q)}$$

where $\mu$ is the *Mobius function* of the poset Q (see Rota [Ro]).

**Proposition 7** If $L$ is a join-distributive lattice with $a$ atoms, then

$$\chi(L,\lambda) = (\lambda - 1)^a \lambda^{r(\hat{1})-a}.$$

*Proof.* Since $\mu(\hat{0},q) = 0$ unless q is a join of atoms ([4], Prop. 2), the only $q \in L$ that contribute to the sum are those in the Boolean algebra $B$ generated by the atoms of $L$. Thus we have

$$\chi(L,\lambda) = \sum_{q \in B} \mu(\hat{0},q) \lambda^{a-r(q)} \lambda^{r(\hat{1})-a}$$

$$= \chi(B,\lambda) \lambda^{r(\hat{1})-a}$$

$$= (\lambda-1)^a \lambda^{r(\hat{1})-a}$$

using the well-known fact (see e.g. [Ro]) that $\chi(B,\lambda) = (\lambda - 1)^a$.□

Using proposition 7, we get the following form for $\chi(L,\lambda)$:

**Corollary 8** Let $s$ be the number of elements of $P$ which are covered by a unique element, and let $m$ be the cardinality of $\overline{M}$ (= the mlb-closure of the maximal elements of $P$). Then

$$\chi(LC(P),\lambda) = (\lambda - 1)^s \lambda^{card(P)-m-s}.$$

*Proof.* By Corollary 6 and Theorem 4 we have $r(\hat{1}) = card(P) - m$ and so to apply Proposition 7 we only need to show that the atoms of $LC(P)$ correspond to elements of $P$ covered by a unique element. Let $H$ be an atom of $LC(P)$ with $A$ the se t of its mlb-closed elements. By Corollary 6, $card(P) - card(A) = r(H) = 1$ so there is exactly one non-closed element $x$ for $H$. Hence $H(x)$ covers $x$, and if any other element $y$ covers $x$, then $H(y) \geq H(x)$ implies $H(y) \neq y$, a contradiction.

Conversely, specifying $x$ to be the only unclosed element does define a closure relation. □

**Example.** Figure 1 shows the atoms of $LC(P)$, the values of $\mu(\hat{0},x)$ and $\chi(LC(P),\lambda)$.

# 3 Modular elements and supersolvability

Using Theorem 4, it is easy to characterize which closure relations correspond to various distinguished classes of elements of $LC(P)$, such as the atoms, coatoms, join-irreducibles, and meet-irreducibles. One interesting class for which this is non-trivial are the *modular* elements of $LC(P)$. Recall ( see Stanley [St]) that an element $H$ in a lattice $L$ is modular if and only if for all $K \leq K' \in L$ we have

$$H \vee K = H \vee K', \text{ and } H \wedge K = H \wedge K' \text{ imply } K = K'.$$

**Theorem 9.** Assume $P$ has a greatest element $\hat{1}$. Then a closure relation $H$ on $P$ is a modular element of $LC(P)$ if and only if $H$ satisfies the following *cover property*:

For all $x, y \in P$, if $y$ covers $x$ and $H(x) \neq x$, then $H(y) = H(x)$.

*Proof.*($\Rightarrow$): Let $H$ be a closure relation on $P$ not having the cover property, i.e. there is some $y$ covering $x$ such that $H(x) \neq x$, but $H(y) \neq H(x)$. Let $K'$ be the closure relation having closed elements $\{\hat{1}, y\} \cup P_{<x}$, and $K$ the closure relation having closed elements $\{\hat{1}, x, y\} \cup P_{<x}$ (note that both of these sets are mlb-closed). Then $K < K'$, and it is easily seen that $H \vee K = H \vee K'$ and $H \wedge K = H \wedge K'$, violating the modularity of $H$.

($\Leftarrow$): We show in general that any closure relation with the cover property is modular. Suppose $H$ satisfies the cover condition, and assume we have a pair of closure relations $K \leq K'$ such that $H \vee K = H \vee K'$ and $H \wedge K = H \wedge K'$. We must show that $K = K'$, so it would suffice to show $K \geq K'$. So assume $K(x) = x$ for a given $x$, and we will prove that $K'(x) = x$. From $H \vee K = H \vee K'$, we may assume that $H(x) \neq x$. Our strategy: We claim that any $y \geq x$ must be comparable to $H(x)$. Next we use this claim to show $K'(x) < H(x)$ by contradiction. Finally, we show $K'(x) = x$.

To prove the claim, suppose $y \not\leq H(x)$. Then $H(y) \neq H(x)$, so considering a maximal chain from $x$ to $y$, there must exist $y''$ covering $y'$ such that $H(y') = H(x)$ but $H(y'') \neq H(x)$. Now by the cover property, we must have $H(x) = H(y') = y'$. Therefore $y \geq y' = H(x)$ and the claim is proved.

Our next goal is to show that $K'(x) < H(x)$. By the claim we know that these two elements must be comparable. So suppose for contradiction that $K'(x) \geq H(x)$. Let $K'_c, H_c$ denote the sets of closed elements of $K', H$ respectively. Since $min(H_c)_{\geq x} = H(x)$ and $min(K'_c)_{\geq x} = K'(x) > H(x)$, any maximal lower bound $z \geq x$ of a subset $A \subseteq H_c \cup K'_c$ must be comparable to $H(x)$ (by the claim), and hence $z \geq H(x)$ by maximality. Thus

$$(H \wedge K')(x) = min(\overline{H_c \cup K'_c})_{\geq x} \geq H(x) > x,$$

contradicting the fact that $(H \wedge K')(x) = (H \wedge K)(x) = x$. Hence $K'(x) < H(x)$.

Lastly, we show $K(x) = x$. Let $k' \in K'_c$ and $h \in H_c$ have a maximal lower bound $z \geq x$. By our claim either $k' \geq H(x)$, in which case $z \geq H(x)$, or $k' < H(x)$, in which case $z = k'$. This shows that

$$(\overline{K'_c \cup H_c})_{<H(x)} = (K'_c)_{<H(x)}.$$

But

$$x = (H \wedge K)(x) = (H \wedge K')(x) \in (\overline{K'_c \cup H_c})_{<H(x)},$$

so $x \in K'_c$, i.e. $K'(x) = x$.□

As a corollary, we have

Figure 2: Constructing an $M$-chain in $LC(P)$.

**Theorem 10** $LC(P)$ is supersolvable, i.e. it contains a maximal chain of modular elements, (See Stanley [St] for alternate definitions and consequences of supersolvability).

*Proof.* We first prove the theorem with the assumption that $P$ has a greatest element $\hat{1}$, and then deduce the theorem in general.

Let $p_1, p_2, \ldots, p_n$ be any linear extension of the partial order on $P$. Let $H_0$ be the greatest element of $LC(P)$, having $\hat{1}$ as is its only closed element, and for $i = 1, 2, \ldots, n$ let $H_i$ be the closure relation having closed elements $\{\hat{1}, p_1, p_2, \ldots, p_i\}$. One can check that this set is in fact mlb-closed, since any maximal lower bounds of subsets of $\{\hat{1}, p_1, p_2, \ldots, p_i\}$ must either be $\hat{1}$ or lie in the ideal $\{p_1, p_2, \ldots, p_i\}$. It is easy to see that each $H_i$ satisfies the cover condition of Theorem 9, and hence is a modular element of $LC(P)$: if $y$ covers $x$ and $H_i(x) \neq x$, then $x \notin \{\hat{1}, p_1, p_2, \ldots, p_n\}$, so either $y = \hat{1}$ or else $y \notin \{\hat{1}, p_1, p_2, \ldots, p_n\}$. In either case, $H_i(x) = \hat{1} = H_i(y)$.

To prove the theorem in general, let $P$ be an arbitrary poset with maximal elements $M$. Adjoin a greatest element $\hat{1}$ to $P$ to obtain the poset $P + \hat{1}$. Let $H$ be the closure relation on $P + \hat{1}$ with closed elements $\{\hat{1}\} \cup \overline{M}$ (an mlb-closed set). Then by Theorem 4, we have that $LC(P)$ is isomorphic to the interval between the identity closure and $H$ in $LC(P + \hat{1})$. Since $LC(P + \hat{1})$ is supersolvable, this interval is also supersolvable ([St], Proposition 3.29(i)), and hence $LC(P)$ is supersolvable.□

**Acknowledgements.** The authors would like to thank Mark Haiman and G.C.-Rota for several helpful discussions and the referee for numerous suggestions and corrections.